%% file: art5.tex
\documentclass[10pt, a4paper, reqno]{amsart}

\usepackage{amsthm, amsmath, amssymb, amsfonts, stmaryrd, latexsym, mathrsfs, enumerate}
\usepackage[english]{babel}
\usepackage{graphicx}
\usepackage[all]{xy}
\usepackage{nicefrac}
\usepackage{color}
\usepackage{mparhack}
\usepackage{algorithmic}
\usepackage{algorithm}

\theoremstyle{plain}
\newtheorem{lemma}{Lemma}
\newtheorem{theo}{Theorem}
\newtheorem{prop}{Proposition}
\newtheorem{remark}{Remark}

\theoremstyle{definition}
\newtheorem{definition}{Definition}

\newcommand{\F}{\mathbb{F}}
\newcommand{\Z}{\mathbb{Z}}
\newcommand{\Q}{\mathbb{Q}}
\newcommand{\R}{\mathbb{R}}

\newcommand{\kp}{\mathfrak{p}}

\newcommand{\sD}{\mathscr{D}}

\newcommand{\kk}{\kappa}

\newcommand{\cO}{\mathcal{O}}
\newcommand{\Gal}{\mathop{\mathrm{Gal}}\nolimits}
\newcommand{\Aut}{\mathop{\mathrm{Aut}}\nolimits}

\newcommand{\sfr}[2]{\nicefrac{#1}{#2}}

\newcounter{enumi_saved}

\makeatletter
\def\imod#1{\allowbreak\mkern10mu({\operator@font mod}\,\,#1)}
\makeatother

\allowdisplaybreaks

\begin{document}

\title[A family of Eisenstein polynomials]{A family of Eisenstein polynomials generating totally
  ramified extensions, identification of extensions and construction of class fields} 
\author{Maurizio Monge}
\email{maurizio.monge@sns.it} 
\address{Scuola Normale Superiore di Pisa - Piazza dei Cavalieri, 7 - 56126 Pisa} 
\date{\today}

\keywords{Eisenstein polynomial, normal form, totally ramified extension, special generator, p-adic
  field, ramification theory, algorithm, local class field theory, Serre mass formula}
\subjclass[2010]{11S15, 11Y40, 11S31}

\begin{abstract}
  Let $K$ be a local field with finite residue field, we define a normal form for Eisenstein
  polynomials depending on the choice of a uniformizer $\pi_K$ and of residue representatives. The
  isomorphism classes of extensions generated by the polynomials in the family exhaust all totally
  ramified extensions, and the multiplicity with which each isomorphism class $L/K$ appears is
  always smaller than the number of conjugates of $L$ over $K$.

  An algorithm to recover the set of all special polynomials generating the extension determined by
  a general Eisenstein polynomial is described. We also give a criterion to quickly establish that a
  polynomial generates a different extension from that generated by a set of special polynomials,
  such criterion does not only depend on the usual distance on the set of Eisenstein polynomials
  considered by Krasner and others.

  We conclude with an algorithm for the construction of the unique special equation determining a
  totally ramified class field in general degree, given a suitable representation of a group of
  norms.
\end{abstract}

\maketitle

\section{Introduction}

In this note we show how it is possible to define a normal form for Eisenstein polynomials, which
can be used for quickly enumerating totally ramified extensions of a local field, for selecting a
special defining polynomial to represent extensions, and for identification of the
extensions. Unluckily it doesn't seem possible to produce easily exactly one special polynomial for
each isomorphism class of extensions, but we show how to obtain a very restricted set of polynomials
generating each extension. The number of special polynomials generating a fixed extension $L/K$ is
always smaller that the number of conjugates of $L$ over $K$, so that each Galois extensions is
generated by exactly one polynomial. In fact, the problem of selecting exactly one generating
polynomial for each isomorphism class appears to be as hard as that of determining the cardinality
of the group of automorphisms of the extension generated by a polynomial.

As shown in \cite{Pauli2001computation}, it is possible to enumerate and identify the extensions
generated by Eisenstein polynomials selecting one polynomial for each neighborhood with respect to a
suitable distance, and applying Panayi root finding algorithm to collect the polynomials which
generate the same extension. The search space can be drastically reduced by just taking into account
Eisenstein polynomials in normal form.

Furthermore, for each Eisenstein polynomial generating an extension $L/K$ of degree $n$ there exists
a quick way to recover all the special polynomials attached to the extension, which does not require
an exhaustive search over the space of all extensions of degree $n$ of $K$, not even a search within
the set of polynomials generating extensions with fixed ramification data.

Indeed, any Eisenstein polynomial can be put into normal form by applying greedily a reduction
algorithm, which however allows some free choices during the reduction. The full set of special
polynomials is obtained as the set of all possible outputs of the reduction algorithm, over all
possible choices.

We exhibit a criterion for establishing \emph{a priori} that an Eisenstein polynomial $f(T)$ may not
be converted to another polynomial $g(T)$ via such a reduction applied greedily, and when $f(T)$ and
$g(T)$ are any two Eisenstein polynomials such that one of them is known to generate a Galois
extension then the criterion can be used to show that $f(T)$ and $g(T)$ generate non-isomorphic
extensions. The criterion takes into account the higher order terms appearing in the $p$-adic
expansion of the coefficients, not just the valuation (or first-order expression) of $f(\pi)-g(\pi)$
for a uniformizer $\pi$ of an extension $L/K$ of degree $n$. The criterion established in
\cite{yoshida2011ultrametric} for totally ramified Galois extensions over $\Q_p$ is also recovered
in a more general context.

In the last section we describe an algorithm which allows to construct the unique special Eisenstein
polynomial generating a totally ramified class field, given a suitable description of a norm
subgroup. In particular, we show that there exists an ordering of the terms appearing in the
$p$-adic expansions of the coefficients which allows to recover all the terms of the special
polynomial, by solving inductively linear equations over the residue field. An algorithm for the
construction of polynomials generating class field was described in \cite{pauli2006constructing} for
cyclic extensions, where an extension of degree $p^m$ is constructed inductively by steps of degree
$p$. In our construction an Eisenstein polynomial generating an arbitrary totally ramified class
field is constructed directly.

\subsection{Acknowledgements}
We would like to thank Philippe Cassou-Nogu\`es, Ilaria Del Corso, Roberto Dvornicich and Boas Erez
for various discussions on this topic, and the Institut de Math\'ematiques de Bordeaux for
hospitality while conceiving this work.

\subsection{Notation} All local fields will be assumed to have finite residue field, even though
most results hold just assuming it to be perfect, and some even in greater generality.

If $F$ is any local field, we denote by $\cO_F$ its ring of integers with unique maximal ideal
$\kp_F$, and by $\kk_F$ the residue field $\cO_F/\kp_F$. We denote by $v_F$ the normalized valuation
and by $e_F$ the absolute ramification index $v_F(p)$ (if $F$ has characteristic $p$ we will be fine
with $e_F=+\infty$).

Let $L/K$ be a separable totally ramified extension of degree $n$, we denote by $\Gamma=\Gamma(L/K)$
the set of all $K$-embeddings of $L$ into a fixed separable algebraic closure $L^{\mathrm{sep}}$ of
$L$. For $\sigma\in{}\Gamma(L/K)$ we define
\[
    i_{L/K}(\sigma) = \min_{x\in{}\cO_L}v_L(\sfr{\sigma(x)}{x}-1),
\]
which is also equal to $v_L(\sfr{\sigma(\pi)}{\pi}-1)$ for a uniformizing element $\pi=\pi_L$
generating the maximal ideal $\kp_L$ of $\cO_L$. For each real $x$ we set
\begin{align*}
   \Gamma_x &= \{\sigma\in\Gamma : i_{L/K}(\sigma)\geq{}x\}, \\
   \Gamma_{x^+} &= \{\sigma\in\Gamma : i_{L/K}(\sigma)>{}x\},
\end{align*}
our definition is equal to that in \cite{Yamamoto1968} and reduces to the classical definition of
ramification subgroups of \cite{fesenko2002local,serre1979local,deligne1984lescorps} for Galois
extensions, and differs by a shift by $1$ from that used in \cite{lubin1981local,Helou1991}.

We say that $x$ is a \emph{ramification break} if $\Gamma_x\supsetneq\Gamma_{x^+}$, and let $K_{x}$
and $K_{x^+}$ respectively be the fields fixed by $\Gamma_x$ and $\Gamma_{x^+}$.

We put
\[
   \phi_{L/K}(x) = \frac{1}{n}\int_0^x \left(\#\Gamma_t\right) dt,
\]
and let $\psi_{L/K}$ to be the inverse by composition. Since we only consider totally ramified
extensions we restrict both $\phi_{L/K}$ and $\psi_{L/K}$ to $\R^{\geq0}$.

If $L/K$ is any totally ramified extension of degree $n$, with $k$ distinct ramification jumps say,
we will usually denote with $t_1<t_2<\dots<t_k$ the ramification jumps and with
$\gamma_0>\gamma_1>\dots>\gamma_k$ the cardinalities of the corresponding ramification subsets, so
that $\gamma_0=\#\Gamma=n$ and $\gamma_i=\#\Gamma_{t_i^+}$ for $1\leq{}i\leq{}k$. The $\gamma_i$ are
all powers of $\pi$, except possibly for $\gamma_0=n$.

If $p^s$ is the biggest power of $p$ dividing $n$, for each $0\leq\ell\leq{}s$ it will also be
convenient putting $\tau_\ell$ to be equal to the smallest real $t$ such that $n\phi_{L/K}(x)$ has
slope $\leq{}p^\ell$ for $x\geq{}t$, it will be equal to either $0$, or some ramification jump
$t_i$. The $\tau_\ell$ are weakly decreasing and exhaust all the lower ramification jumps $t_i$, and
one jump $t_i>0$ is repeated $r$ times if $(\Gamma_{t_i}:\Gamma_{t_i^+})$ is equal to $p^r$, so each
ramification jump is taken ``with its multiplicity'' in a suitable sense. It will also be convenient
defining
\[
     \xi_\ell = n\phi_{L/K}(\tau_\ell)-p^\ell\tau_\ell,\qquad
     \sigma_\ell = n\phi_{L/K}(\tau_\ell),
\]
for each $0\leq\ell\leq{}s$. Up to a factor $n$ the $\sigma_\ell$ are the upper ramification jumps
of the extension.

\section{Reduction algorithm and the family of reduced polynomials}

Let $f(T)=T^n + f_{n-1}T^{n-1} \dots f_1T + f_0$ be a monic Eisenstein polynomial of degree $n$, let
$\pi$ be a root in a fixed algebraic closure $K^\mathrm{alg}$ and put $L=K(\pi)$. Then clearly
$f(T)$ is the minimal polynomial of $\pi$ which is a uniformizing element of the extension
determined by $f(T)$, and we are interested in understanding how the coefficients of the minimal
polynomial of a uniformizer change when $\pi$ is replaced by another uniformizer
$\rho=\pi+\theta\pi^{m+1}+\cO(\pi^{m+2})$, for some unit $\theta\in{}U_K$ and integer
$m\geq1$. Since the computation which follows only depends on $\theta$ at the first order, $\theta$
may be taken to be a multiplicative representative.

Let us consider the ramification polynomial $\Phi(T)=\pi^{-n}f(\pi{}T+\pi)$, its Newton polygon is
fully described by the lower ramification breaks. For $\alpha\in{}\cO_K$ we can compute a lower
bound for the valuation of $\Phi(\alpha)$ as function of $v_L(\alpha)$ starting from the Newton
polygon of $\Phi(T)$. The construction produces naturally the \emph{Newton copolygon}, which is
essentially the dual convex body of the Newton polygon, and is connected to the Hasse-Herbrand
transition function as already observed in \cite{lubin1981local,li1997p}; in such references
$f(T+\pi)$ was used instead so the function obtained was slightly different from the classical
Hasse-Herbrand defined in \cite{fesenko2002local,serre1979local}.

Indeed, the Newton polygon of the polynomial $\Phi(\pi^mT)$ resulting by the substitution
$T\rightarrow{}\pi^mT$ can obtained from the polygon of $\Phi(T)$ moving the points with abscissa
$x$ up by $\frac{m}{n}x$. In other words, if $N:[1,n]\rightarrow\R$ is the real function
describing the polygon of $\Phi(T)$, the polygon of $\Phi(\pi^mT)$ is described by
$N(x)+\frac{m}{n}x$.

The function $N(x)$ is convex and piecewise linear, and by the well known properties of Newton
polygons the slopes are $-\sfr{t_k}{n},\dots,-\sfr{t_1}{n}$ where $t_1<t_2<\dots<t_k$ are the lower
ramification breaks of the extension generated by a root, and it has slope $-\sfr{t_i}{n}$ in the
interval $[\gamma_i,\gamma_{i-1}]$ where $\gamma_0>\gamma_1>\dots>\gamma_k$ are the cardinalities of
the corresponding ramification subsets. We put $t_0=+\infty$, $t_{k+1}=-\infty$ for convenience.

\begin{center}
\scalebox{0.9}{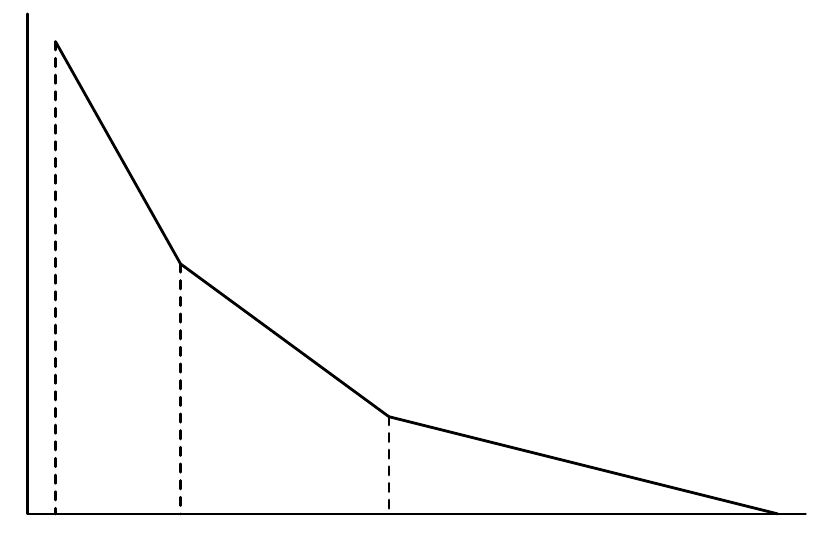}
\end{center}

Let's consider the minimum achieved by the function $N(x)+\frac{m}{n}x$ in the interval $[1,n]$, as
a function of the real parameter $m$. It is again a piecewise linear function with slope
$\sfr{\gamma_i}{n}$ for $t_i\leq{}m\leq{}t_{i+1}$, and we obtain that this minimum value function is
exactly the Hasse-Herbrand function $\phi_{L/K}(m)$. Hence this is the smallest valuation (with
respect to $K$) of the coefficients of $\Phi(\pi^mT)$, and $\pi^{-n\phi_{L/K}(m)}$ is the exact
power of $\pi$ such that $\pi^{-n\phi_{L/K}(m)}\Phi(\pi^mT)$ is in $\cO_L[T]$ and non trivial modulo
$\kp_L$.

Let's define the valuation of a polynomial to be the smallest valuation of the coefficients, we can
resume what proved in the following

\begin{prop}
\label{sec1_prop1}
  Let $f(T)$ be an Eisenstein polynomial, $\pi$ a root, $L=K(\pi)$ and
  $\Phi(T)=\pi^{-n}f(\pi{}T+\pi)$ its ramification polynomial. Then
\[
   v_L(\Phi(\pi^mT)) = n\phi_{L/K}(m).
\]
\end{prop}

It will also be convenient to deduce an expression for the values $N(p^\ell)$ for each $\ell\geq0$
such that $p^\ell\mid{}n$. Starting from $p^\ell$ the function $N(x)$ has slope
$-\sfr{\tau_\ell}{n}$, so $N(x)+\frac{\tau_\ell}{n}x$ has infimum equal to $\phi_{L/K}(\tau_\ell)$,
which is achieved for $x=p^\ell$ and is also equal to $N(p^\ell)+\frac{\tau_\ell}{n}p^\ell$, so we obtain
\[
      N(p^\ell) = \phi_{L/K}(\tau_\ell) - \frac{\tau_\ell}{n}p^\ell=\frac{\xi_\ell}{n}.
\]

\begin{lemma}
\label{sec1_lem1}
For each $\ell\geq0$ such that $p^\ell\|n$ we have $N(p^\ell)=\sfr{\xi_\ell}{n}$.
\end{lemma}

We will also prove another Lemma which we will require later. If $p^\ell$ is the abscissa of a
vertex of the Newton polygon we have that the terms contributing to the coefficient of $T^{p^\ell}$
in the ramification polynomial give to the coefficient of $T^{p^j}$ contributions having
$K$-valuation at most $e_K(\ell-j)$ bigger, for $j<\ell$. In other words we have that
$N(p^j)<e_K(\ell-j)+N(p^\ell)$, for each $j\leq\ell$. Considering the last vertex of one side of the
Netwon polygon, and since for each $\ell$ the slope is equal to $-\sfr{\tau_\ell}{n}$ in the
interval $[p^\ell,p^{\ell+1}]$ and $N(p^\ell)<e_K+N(p^{\ell+1})$, we have that then $\tau_\ell$ has
to be at most $n\frac{e_K}{(p^{\ell+1}-p^\ell)}=\sfr{e_L}{(p^{\ell+1}-p^\ell)}$. Hence we have

\begin{lemma}
\label{sec1_lem2}
We have
\[
    \xi_j \leq e_L(\ell-j) + \xi_\ell
\]
for each $j<\ell$, and furthermore
\[
   \tau_\ell \leq \sfr{e_L}{(p^{\ell+1}-p^\ell)}
\]
for each $\ell$.
\end{lemma}

We now study the points $(j,v_K(\Phi_j))$ coming from a monomial $\Phi_jT^j$ that may lie on the
boundary of the Newton polygon of $\Phi(T)=\sum_{i=0}^n\Phi_iT^i$. We claim that either their
ordinate $j$ is a power of $p$, either $\gamma_1\mid{}j$, and the latter is only possible when
$\gamma_0=n$ is not a power of $p$, so that $\gamma_1$ is the biggest power of $p$ dividing $n$, and
the polygon of $\Phi(T)$ has slope $0$ in the interval $[\gamma_1,\gamma_0]$.

Indeed, for each $\ell$ we have
\[
  \Phi_{p^\ell} = \sum_{i=p^\ell}^n\binom{i}{p^\ell}f_i\pi^{i-n},
\]
and since the summands have different valuations modulo $n$ the valuation of $\Phi_{p^\ell}$ has to
be equal to the minimal valuation of such terms. For any integer $r$ the terms
$\binom{i}{r}f_i\pi^{i-n}$ contributing to $\Phi_{r}$ have valuation which is at least as big as the
valuation of $\binom{i}{p^\ell}f_i\pi^{i-n}$ when $p^{\ell+1}\nmid{}r$, and strictly bigger if
$p^\ell\nmid{}r$, so when $p^\ell\|r$ we have that $v_L(\Phi_{r})\geq{}v_L(\Phi_{p^\ell})$ and
$(r,v_K(\Phi_r))$ cannot be on the boundary of the polygon unless possibly when the segment
containing $p^\ell$ has horizontal slope, $p^\ell=\gamma_1$ and $p^\ell\mid{}r$.

For integer $m\geq0$ let's consider the polynomials
\[
   S_m(T) = \overline{\pi^{-n\phi_{L/K}(m)}\Phi(\pi^mT)}.
\]
If $m\geq1$, or $n$ is a power of $p$, then $S_m(T)$ is of the form
\[
   S_m(T) = \sum_{i=a}^b c_iT^{p^i}
\]
for some coefficients $c_i$, where $p^a=p^b=\gamma_i$ when $m$ is not a ramification break and
$m\in(t_{i+1},t_{i})$ say, while $\gamma_i=p^a$ and $\gamma_{i-1}=p^b$ if $m=t_i$ for some $i$. In
particular they are additive polynomials.

On the other hand if $m=0$ and $n$ is not a power of $p$ (and hence $L$ has a non-trivial tamely
ramified subextension) the terms appearing in $S_0(T)=\overline{\Phi(T)}$ are all coming from the
leading monomial $T^n$, so that putting $n'=\sfr{n}{\gamma_1}$ we have
\begin{align*}
  S_0(T) &= \sum_{j=1}^{n'}\binom{n}{\gamma_1j}T^{\gamma_1j} \\
   &=  \sum_{j=1}^{n'}\binom{n'}{j}T^{\gamma_1j} \\
   &= (1+T^{\gamma_1})^{n'} - 1.
\end{align*}
We collect these facts in the following proposition.

\begin{prop}
  If $m\geq1$ the polynomial $S_m(T)$ is an additive polynomial, which is composed by more than one
  monomial if and only if $m$ is a lower ramification break. For $m=0$ we have
  $S_0(T)=(1+T^{p^s})^{n'}-1$, where $n=p^sn'$ and $(p,n)=1$.
\end{prop}

When the context is clear, we will abuse of notation and also denote by $S_m$ the induced map
$\theta\mapsto{}S_m(\theta)$ over the residue field or an extension thereof.

\subsection{Change induced on the coefficients by a substitution}
We study now the effect of replacing the minimal polynomial $f(T)$ of $\pi$ with the minimal monic
polynomial $g(T)$ of a different uniformizer $\rho$.

Let's take $\rho=\pi+\theta\pi^{m+1}+\cO(\pi^{m+2})$, we will identify the term $(f_i-g_i)\rho^i$
which has minimal valuation for general $\theta$, and which gives information about the most
significant change induced on the coefficients $f_i\rightarrow{}g_i$ as consequence of the
substitution $\pi\rightarrow\rho$.

The non-zero terms $(f_i-g_i)\rho^i$ have valuations with different remainders modulo $n$, and
furthermore we have
\begin{align*}
  \sum_{i=0}^{n-1} (f_i-g_i)\rho^i &= f(\rho) - g(\rho) \\
      &= f(\rho) = \pi^n\Phi(\theta\pi^m+\cO(\pi^{m+1})),
\end{align*}
considering the definition of $\rho$. If $m\geq1$, being $\rho\equiv\pi\imod{\kp_K^2}$ we obtain the
following Lemma, after dividing by $\pi^{n(\phi_{L/K}(m)+1)}$ and reducing the expression modulo
$\kp_K$.

\begin{lemma}
\label{sec1_lem3}
If $m\geq1$ and $g(T)$ is the minimal monic polynomial of an element of the form
$\rho=\pi+\theta\pi^{m+1}+\cO(\pi^{m+2})$ we have
\[
   \overline{(f(\pi)-g(\pi))\cdot \pi^{-n(\phi_{L/K}(m)+1)}} = S_m(\theta).
\]
\end{lemma}

Since $n\mid{}v_L(f_i-g_i)$ for each $i$, the unique term $(f_i-g_i)\pi^i$ of $f(\pi)-g(\pi)$ which
may be contributing to the left hand side is for $i$ satisfying
\[
   i \equiv n(\phi_{L/K}(m)+1) \imod{n},
\]
so $i$ is uniquely determined being $0\leq{}i<n$. We observe that if $m\geq1$ is not a lower
ramification break then $S_m$ is surjective being $\kk_K$ finite and hence perfect, while if $m=t_i$
for some $i$ then it may not be surjective, when the additive polynomial $S_{t_i}(T)$ has a root
over the residue field $\kk_K$.

Assume $t_i$ to be an integer, we will later show that the polynomial $S_{t_i}(T)$ only depends on
the field extension $L/K$ and on the class of $\pi\mod{}\kp^2$, as a consequence of a stronger
result, Theorem \ref{theo2}, which is proved independently. For the moment we can give a definition
of reduced polynomial without assuming this invariance, even though the definition will be less
manageable from a practical point of view.

Let $I_m$ be the image of $S_m$, and also its preimage in $\cO_L$ when the context is clear. Lemma
\ref{sec1_lem3} says that passing to the minimal polynomial of an element of the form
$\pi+\theta\pi^{m+1}+\cO(\pi^{m+2}))$, if $n(\phi_{L/K}(m)+1)=jn+i$ with $0\leq{}i\leq{}n$, we can
change the corresponding term $f_i$ by an element of $\pi^{nj}I_m$ while all other terms $f_r\pi^r$
are unchanged modulo $\pi^{jn+i+1}$. Since the polynomials $S_r$ for $r\leq{}m$ are certainly
unchanged too, this observation motivates the following definition.

\begin{definition}
\label{def_red}
Let $f(x)$ be an Eisenstein polynomial, and assume each coefficient $f_i$ to have an expansion
\[
     f_i = \sum_{j\geq1}^\infty f_{i,j} \pi_K^j
\]
with $f_{i,j}\in{}R$ for a fixed set of residue representatives $R$, and where $\pi_K$ is a fixed
uniformizer of $K$, and let $\bar\eta_f=\overline{\sfr{-f_0}{\pi_K}}$. Assume the choice of a set
$A_0\subset{}\kk_K^\times$ of representatives of $\kk_K^\times/(\kk_{K}^\times)^n$, and for each
additive polynomial $S_m(T)$ for $m\geq1$ a set of elements $A_m\subset{}\kk_K$ which are
representatives of the cokernel of the map $\theta\mapsto{}\eta_f^jS(\theta)$, where
$j=[\phi_{L/K}(m)+1]$.

We say that $f(x)$ is \emph{reduced} (with respect to the choice of the $A_i$) when we have
\begin{enumerate}
\item\label{c1} $\bar\eta_f=\overline{-\sfr{f_0}{\pi_K}}$ is in $A_0$,
\item\label{c2} for each $m\geq1$, if $n(\phi_{L/K}(m)+1)=jn+i$ for positive integers $i,j$ with
  $i<n$, then we have $\overline{f_{i,j}}\in{}A_m$
\end{enumerate}
We say that $f(x)$ is \emph{reduced up to the level} $r$ when condition \ref{c1} is satisfied, and
condition \ref{c2} holds for all $m\leq{}r$.
\end{definition}

If $f(T)$ is any Eisenstein polynomial, it's easy to see that the polynomial
$\theta^{n}f(\theta^{-1}T)$ satisfies condition \ref{c1} for some suitable $\theta$. A polynomial
reduced up to level $0$ can be obtained by the following algorithm.
\begin{algorithm}[H]
\caption{Reduction (step $0$)}
\label{reduction_m}
\begin{algorithmic}
\STATE $\bar\alpha\gets\overline{\sfr{-f_0}{\pi_K}}$,
\STATE $\bar\beta\gets{}Representative(\bar\alpha, (\kk_K^\times)^n)$,
\STATE $\bar\theta\gets{}Solve(T^n=\bar\beta/\bar\alpha)$,
\STATE $\theta\gets{}Lift(\bar\theta)$,
\RETURN $\theta^{n}f(\theta^{-1}T)$.
\end{algorithmic}
\end{algorithm}

If $i,j$ are as above, we have shown above $f_{i,j}\pi_K^j$ can be changed by any element in
$\pi^{nj}I_m$ modulo $\pi^{nj+1}$, so $\overline{f_{i,j}}$ can be changed by any element of
$\overline{(\sfr{\pi^n}{\pi_K})^j}I_m$. Since $\pi^n=-f_0+\cO(\pi^{n+1})$ we have that
$\overline{(\sfr{\pi^n}{\pi_K})}=\bar\eta_f$, so $\overline{f_{i,j}}$ is changed by an element of
$\bar\eta_f^j{}I_m$, and note that $\eta_f$ is unchanged when passing to the minimal polynomial of an
uniformizer of the form $\pi+\theta\pi^{m+1}+\cO(\pi^{m+2}))$, for a suitable $\theta$.

In particular, if $f(x)$ is reduced up to the level $m-1$ we can obtain a polynomial reduced up to
the level $m$ via the following reduction step.
\begin{algorithm}[H]
\caption{Reduction (step $m$)}
\label{reduction_m}
\begin{algorithmic}
\STATE $j\gets\lfloor\phi_{L/K}(m)+1\rfloor$,
\STATE $i\gets{}n\cdot\{\phi_{L/K}(m)+1\}$,
\STATE $\bar\alpha\gets{}\overline{f_{i,j}}$,
\STATE $\bar\beta\gets{}Representative(\bar\alpha, \mathrm{image}(\bar\eta_f^jS_m))$,
\STATE $\bar\theta\gets{}Solve(\bar\eta_f^jS_m(T) = \bar\alpha-\bar\beta)$,
\STATE $\theta\gets{}Lift(\bar\theta)$,
\STATE $F(T)\gets{}T+\theta{}T^{m+1}+\{\text{any terms of degree}\geq{}m+2\}$,
\RETURN $Resultant_U(f(U), T-F(U))$.
\end{algorithmic}
\end{algorithm}

Indeed, if $g(T)$ is the returned polynomial we have
\[
   \overline{(f_{i,j}-g_{i,j})\pi_K^j\pi^{i-n(\phi_{L/K}(m)+1)}} =
   \overline{(f_{i,j}-g_{i,j})}\bar\eta_f^{-j} = S_m(\bar\theta),
\]
and consequently $\overline{g_{i,j}}=\bar\beta\in{}A_m$.  Since we allow any higher order term in
the choice of $F(T)=T+\theta{}T^{m+1}+\dots$, we anticipate that for a suitable $F(T)$ it will not
be necessary to compute the resultant appearing in the algorithm as the determinant of a big matrix
with coefficients in $K[T]$, see Remark \ref{rmk_resultant}.

\begin{remark}
  If $m$ is bigger than the biggest lower ramification break $t_k$, then $S_m(x)$ is surjective, and
  the function $n(\phi_{L/K}(m)+1)$ assumes as possible values all integers
  $>n(\phi_{L/K}(t_k)+1)$. Consequently we can arbitrarily change all the representatives $f_{i,j}$
  whenever
\begin{equation} \label{krasnquant}
   v_L(\pi_K^j\pi^i)=nj+i > n(\phi_{L/K}(t_k)+1),
\end{equation}
without affecting the generated extension, turning them all to $0$ for instance. In this way we
recover the well known quantitative criterion on the distance of two Eisenstein polynomials ensuring
that they generate the same extension, as considered in
\cite{Krasner1962,Pauli2001computation,yoshida2011ultrametric}.
\end{remark}

\subsection{Characterizing reduced polynomials}
We start with a few remarks about Definition \ref{def_red}. Since we allow a different
choice of the representing sets $A_m$ for each $m$, where the $0$ element of the image of the map is
not even requested to be represented by $0$, we have that each Eisenstein polynomial is reduced for
a suitable choice of the $A_m$. While this choice is very far from what would be recommended in a
computer algebra system it will be useful to be able to consider each Eisenstein polynomial as being
already reduced.

On the other hand on a computer algebra system we can expect to have a more or less canonical way
for selecting representing elements of a quotient, and selecting $0$ as representative of the zero
element in the quotient. Under this hypothesis we clarify here how a reduced polynomial looks
like. In particular we will see that, for each $\ell\geq0$ such that $p^\ell$ divides $n$, the
possible valuations of the terms $f_{i,j}\pi_K^{j}\pi^i$ such $p^\ell\|i$ belong to one fixed
interval, \emph{with some exception}.

Fix $\ell$ and let us consider the terms $f_{i,j}\pi_K^jT^i$ with $p^\ell\|i$, we deduce a lower
bound for the value of $nj+i$ from the shape of the Newton polygon of the ramification
polynomial. Indeed, the contribution to the coefficient of $T^{p^\ell}$ in $\Phi(T)$ is
\[
    \pi^{-n}f_{i,j}\pi_K^j\pi^i\binom{i}{p^\ell},
\]
and since the contributions coming from different monomials of $f(T)$ have different valuations
modulo $n$ then their smallest valuation should be at least $nN(p^\ell)=\xi_\ell$. In the same way
we obtain that any term $f_{i,j}\pi_K^jT^i$ with $p^\ell\|i$ and $nj-n+i\geq\xi_\ell$ is compatible
with the ramification data, and when $p^\ell\|\xi_\ell$ and $p^\ell$ is the abscissa of a vertex of
the ramification polygon then there should be a term $f_{i,j}\pi_K^jT^i$ such that the valuation
$nj-n+i$ of the contributed term is exactly $\xi_\ell$, this case corresponds to a vertex of the
Netwon polygon and hence the minimum is reached.

We will show now that starting all the terms $f_{i,j}\pi_K^jT^i$ with $p^\ell\|i$ and $nj+i$ big
enough are turned to $0$ by the reduction algorithm, with a few exceptions.  Indeed, we claim that
the integers which are multiple of $p^\ell$ and $>n\phi(\tau_\ell)=\sigma_\ell$ are all of the form
$n\phi_{L/K}(m)$ for some $m>t_\ell$ (note that $\sigma_\ell$ may not be a multiple of $p^\ell$
itself, we are considering non-Galois extensions and $t_r$ and $\phi(t_r)$ may not be integers).

To show the claim we work by induction on the number of ramification breaks. If
$p^\ell<\gamma_{k-1}$ then $\tau_\ell=t_k$, and $nN(1)$ is certainly an integer being equal to
$v_L(\sD_{L/K})$, and $n\phi_{L/K}(m)$ for integer $m$ assumes as values all integers which are
$>n\phi_{L/K}(t_k)$, being $n\phi_{L/K}(x)$ equal to $nN(1)+x$ for integer $m>t_k$. Assume instead
$p^\ell\geq{}\gamma_{k-1}$, then by induction $\frac{n}{\gamma_{k-1}}\phi_{K_{t_{k-1}^+}/K}(m)$
takes as values any multiple of $p^\ell/\gamma_{k-1}$ which is bigger than
$\frac{n}{\gamma_{k-1}}\phi_{K_{t_{k-1}^+}/K}(\tau_\ell)$ for integer $m>\tau_\ell$. So
$n\phi_{K_{t_{k-1}^+}/K}(x)$ satisfies the required property with respect to $p^\ell$, and so does
$n\phi_{L/K}(x)$ which is obtained as the minimum of $n\phi_{K_{t_{k-1}^+}/K}(x)$ and $nN(1)+x$.

Consequently we have from the claim that all the terms $f_{i,j}\pi_K^j$ with $p^\ell\|i$ and
$nj-n+i\geq{}\sigma_\ell$ can be forced to satisfy $f_{i,j}=0$, except possibly when $nj-n+i$ is
itself equal to $\sigma_r$ for some $r\leq\ell$, in this case we can only force $f_{i,j}$ to be a
suitable representative depending on the image of the polynomial $S_{\tau_r}(T)$, which may not be
surjective as a function over $\kk_K$.

In the case of three breaks we have the following figure representing the values $nj-n+i$ of the
terms of a reduced polynomial.

\begin{center}
\scalebox{0.9}{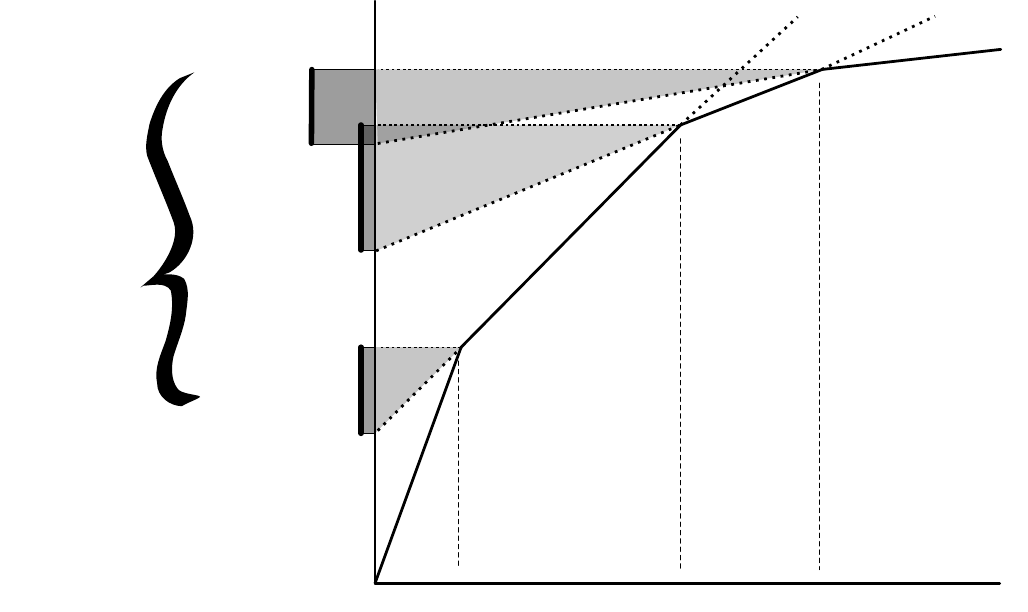}
\end{center}

We state the above results in the following Proposition.

\begin{prop}
\label{prop_range}
Let $f(x)$ be a reduced Eisenstein polynomial, and assume each coefficient $f_i$ to have an expansion
\[
    f_i = \sum_{j\geq1}^\infty f_{i,j} \pi_K^j.
\]
Assume $p^\ell\|i$, then $f_{i,j}$ is non-zero only when
\[
    \xi_\ell \leq nj-n+i < \sigma_\ell,
\]
or when $nj-n+i$ is equal to some $\sigma_r$ and the corresponding additive polynomial
$S_{\tau_r}(T)$ has a root in $\kk_K$.
\end{prop}

In other words we have that starting from a certain points all terms $f_{i,j}\pi_K^j$ for
$p^\ell\|i$ can all be simplified to $0$, except at upper ramification breaks. We will later
see how this phenomenon can be interpreted in terms of local class field theory for abelian
extensions, or in connection with Serre mass formula \cite{serre1978formule} in some simple
particular case.

\subsection{Representation of automorphism as power series}
Applying such substitutions for increasing $m$ we are taking into account all transformations
$F(\pi)$ of $\pi$ by a power series without constant coefficient
$F(T)=\theta_1{}T+\theta_2T^2+\dots$ which may provide an element whose minimal polynomial is
reduced, because any such power series can be written as a composition of polynomials of the form
$T(1+\theta{}T^m)$. Applying the above reduction step for increasing $m$, when $m$ is not equal to a
ramification break $t_i$ we have a unique possible choice for the class $\bar\theta$ of $\theta$ in
the substitution $\pi\rightarrow\pi(1+\theta\pi^m)$. When $m=t_i$ for some $i$, the choice for
$\bar\theta$ is defined up to an element which is a root of $S_{t_i}(T)$, and taking into account
representatives $\theta$ for all possible choices for $\bar\theta$ we can track all possible
outputs. We can run this algorithm starting from the set $\{f(T)\}$ and replacing each polynomial
with the set of all possible outputs, which may not be unique at the ramification breaks $t_i$, and
do so up to the level $t_k$. After this last step we obtain reduced polynomials turning to $0$ all
the $f_{i,j}$ for $i,j$ such that $nj+i>n(\phi_{L/K}(t_k)+1)$.

\begin{algorithm}
\caption{All reduced polynomials}
\label{reduced_multiset}
\begin{algorithmic}
\STATE $\{t_1,\dots,t_k\}\gets{}LowerRamificationBreaks(f(T))$
\STATE $A\gets\{f(T)\}$
\FOR{$m = 0 \to t_k$}
  \STATE $B\gets\emptyset$
  \FOR{$g(T)\in{}A$}
    \STATE $B\gets{}B \cup AllReductions(g(T),m)$
  \ENDFOR
  \STATE $A\gets{}B$
\ENDFOR
\STATE $a\gets[\phi_{L/K}(t_k)+1]$
\STATE $b\gets{}n\cdot\{\phi_{L/K}(t_k)+1\}$
\RETURN $A\mod (\pi^{a+1},\pi^aT^b)$
\end{algorithmic}
\end{algorithm}

Since some outputs may be repeated we end with a multiset of reduced polynomials. Clearly different
power series $F(T)$, $G(T)$ may give the same value $F(\pi)=G(\pi)$ when evaluated in $\pi$, but we
will show that we took into account all the different \emph{values} $F(\pi)\in{}L$ such that the
minimal polynomial of $F(\pi)$ is reduced.

Indeed, in step $0$ we considered all possible values for $F(\pi)$ modulo $\kp_L^2$, and assume by
induction that all the $F(\pi)$ taken in to account up to step $m-1$ cover all possible values
modulo $\kp_L^{m+1}$. The values $F(\pi)+\theta{}F(\pi)^{m+1}+\cO(\pi^{m+2})$ covered in step $m$,
for all admissible representatives $\theta$, will provide all possible values modulo $\kp_L^{m+2}$.

Let $\rho_i(\kk_K)$ be equal to the cardinality of $\kk_K^\times/(\kk_K^\times)^n$ if $t_i=0$, and
to be the number of roots of $S_{t_i}(x)$ contained in $\kk_K$ if $t_i>0$. The cardinality of the
multiset of polynomials obtained as output of the algorithm can be computed counting for each $m$
the number of possible choices which is indeed equal to $\rho_m(\kk_K)$, and the total cardinality
is equal to the product of the $\rho_i(\kk_K)$ over all $i$ such that $t_i$ is an integer, that is
\[
   B_{L/K} = \prod_{\substack{1\leq{}i\leq{}k\\t_i\in{}\Z}} \rho_i(\kk_K).
\]
We give now an interpretation of the $\rho_i(\kk_K)$ as the number of automorphism of intermediate
extensions. Indeed, if $t_i=0$ then $\rho_i(\kk_K)$ counts the number of $n$-th roots of the unity
in $\kk_K$, or equivalently of $n'$-th roots if $n=p^sn'$ with $(n',p)=1$, which is also the number
of automorphisms of a tame extension of degree $n'$ of $K$, like $K_{0^+}/K$ is.

For $t_i>0$ let's consider the intermediate extension $K_{t_i^+}/K_{t_i}$: if $g(T)$ is the minimal
polynomial of $\pi$ over $K_{t_i}$ (which is a factor of $f(T)$) then
\[
  \overline{\pi^{-n(\phi_{i-1}(t_i)+1)}g(\pi^{t_i+1}T+\pi)} = S_{t_i}(T),
\]
and consequently representatives $\theta$ of the roots of $S_{t_i}(T)$ are exactly those such that
\[
    \sigma(\pi)/\pi=1+\theta\pi^{t_i}+o(\pi^{t_i})
\]
for some $K_{t_i}$-automorphism $\sigma\in\Aut(L^{\mathrm{sep}}/K_{t_i})$. Now after extending the
elements of $\Gamma_{t_i^+}$ to the normal closure we have
$\Gamma_{t_i^+}(\sigma|_{L})=\sigma\Gamma_{t_i^+}$, this is immediate considering $\Gamma_{t_i^+}$
as the image of elements of a ramification (normal) subgroup of a bigger Galois extension containing
$L$. Consequently averaging over $\Gamma_{t_i^+}$ we obtain
\[
   \sigma(\pi_{K_{t_i^+}})/\pi_{K_{t_i^+}} = 1 + \theta^{\gamma_i}\pi_{K_{t_i^+}}^{t_i} + o(\pi_{K_{t_i^+}}^{t_i}),
\]
where $\pi_{K_{t_i^+}}=N_{L/K_{t_i^+}}(\pi)$. The equality holds because $t_i$ is smaller than the
all ramification numbers of the extension $L/K_{t_i^+}$, by keeping into account the properties of
the norm map $N_{L/K_{t_i^+}}$ (see \cite[Chap. 3, \S1, Prop. 1.5]{fesenko2002local}).

If $\sigma(\pi_{K_{t_i^+}})\in{}K_{t_i^+}$ then $\bar\theta$ is in $\kk_K$, and on the other hand if
$\bar\theta\in\kk_K$ then $\sigma(\pi_{K_{t_i^+}})$ can be approximated better than any other
conjugate of $\pi_{K_{t_i^+}}$ having $K_{t_i^+}/K_{t_i}$ only one ramification break, and
consequently $\sigma(\pi_{K_{t_i^+}})\in{}K_{t_i^+}$ by Krasner Lemma. In other words we have one
root of $S_{t_i}(T)$ in $\kk_K$ for each conjugate of $\pi_{K_{t_i^+}}$ contained in $K_{t_i^+}$,
and $\rho_i(\kk_K)=\#\Aut(K_{t_i^+}/K_{t_i})$.

So we have that $B_{L/K}$ is an invariant of the extension $L/K$. Considering the subgroups
$\Aut(L/K_{t_i})$ of $\Aut(L/K)$ and the corresponding quotients as subgroups of
$\Aut(K_{t_i^+}/K_{t_i})$, we observe that $B_{L/K}$ provides a ``naive'' upper bound to the cardinality
of $\Aut(L/K)$, but which is in general tighter than the full degree $[L:K]$.

Let $f_1(x),\dots,f_r(x)$ be all the reduced polynomials obtained applying the above algorithm.  The
number of times we obtain the same polynomial $f_i(x)$ is equal to the number of distinct $F_j(\pi)$
such that $f_1(F_j(\pi))=0$, and is consequently equal to the number of roots of $f_i(x)$ contained
in $L$, in another words to the cardinality of $\Aut(L/K)$.

\begin{theo}
\label{theo1}
  Each extension $L/K$ is generated by a reduced polynomial, and the number of reduced polynomials
  generating a fixed extension $L/K$ is
  \[
   B_{L/K}/\#\Aut(L/K).
  \] 
  If $f(x)$ is an Eisenstein polynomial such that a root generates an extension isomorphic to $L$,
  then the reduction algorithm outputs a multiset of cardinality $B_{L/K}$ formed by the reduced
  polynomials, each having multiplicity $\#\Aut(L/K)$.
\end{theo}

We remark that if $F(T)$ is a power series such that $F(\pi)$ is a conjugate of $\pi$, the algorithm
giving the set of special polynomials can collect all the $\theta$ used in the substitutions
$\pi\rightarrow\pi+\theta\pi^{m+1}$ to produce an expression of
\[
    F(T) \mod (f(T),T^{t_k+1}),
\]
which can be used to realize the group $\Aut(L/K)$ as group of truncated power series under
composition, we omit the details of the construction.

Note that there is a unique reduced representative for Eisenstein polynomials generating Galois
extensions, while in general we have a set of polynomials which is equal to the ratio of the
``naive'' bound on the number of automorphisms to the real number of automorphisms. We remark that
extinguishing the redundancy from the above family of reduced polynomials seems to be at least as
hard as computing the cardinality of the automorphism group. This can probably be done in a few
particular cases, like for polynomials of degree $p^2$ over an unramifed extension of $\Q_p$, but
depends indeed on a criterion to detect which extensions are Galois and to establish the cardinality
of the group of automorphisms.

\subsection{Comparison with Amano polynomials and Serre mass formula}

We provide here some qualitative observation, without being completely rigorous. First, if the
degree $n$ is prime with $p$ it's easy to say what are reduced polynomials, and they are all of the
form $T^n+\theta\pi_K$ for some representative $\theta\in{}R$ such that $\bar\theta\in{}A_0$, where
$A_0$ is the chosen set of representatives of $K^{\times}/(K^{\times})^n$.

When $n=p$, Amano defined in \cite{Amano1971} a set of special generating polynomials composed by
trinomials. The equations considered here turn out to look much more complicated ``visually''
because they are no longer trinomials, but the number of parameters is clearly the same, and
nevertheless Amano polynomials do not seem to be easily generalizable to higher degree.

Reduced polynomials of degree $p$ are of the form
\[
    T^p + \sum_{i=1}^p \left( \sum_{\substack{pj+i\geq(p-1)t+p\\pj+i<pt+p}} f_{i,j}\pi_K^j \right)\cdot T^i 
     + \pi_K \left(+ f_{0,{t+1}}\pi_K^{t+1}\right),
\]
for some ramification jump $t$ such that either $t=\sfr{pe_K}{p-1}$, either $t$ is
$<\sfr{pe_K}{p-1}$ and $(p-1)t$ is an integer prime with $p$. Furthermore the term
$f_{0,{t+1}}\pi_K^{t+1}$ is present only when $t$ is an integer and the additive polynomial $S_t(T)$
has a root in $\kk_K$, which is precisely the case of the extensions being Galois.

We remark that given an extension $L/K$ of degree $p$, then in the Galois cyclic case the
uniformizer $\pi_K$ may not be a norm by class field theory, so in general an additional term is
indeed required. On the other hand if $L/K$ is not Galois then $\pi_K$ is always in
$N_{L/K}(L^\times)$.

We give one last interpretation of this fact, under the light of the proof of Serre ``mass
formula''. Considering the map
\[
  \left\{ \parbox{1.2in}{uniformizers of extensions of degree $p$} \right\}
    \stackrel{\text{``minimal polynomial''}}{\longrightarrow}
  \left\{ \parbox{1.2in}{Eisenstein polynomials of degree $p$} \right\}
\]
we have a $p$-to-$1$ correspondence between measure spaces, whose scaling factor turns out to be
determined by the discriminant of the extensions as proven in \cite{serre1978formule}. Let's
restrict the map to the uniformizers of a fixed extension $L/K$ in the algebraic closure, then
either the extension is Galois and the map is still $p$-to-$1$, either the extension is not Galois
and the map becomes $1$-to-$1$, but in this case the image has bigger measure.

In other words, for fixed degree and restricting to extensions with a fixed discriminant, the
smaller is the space of polynomials generating one fixed isomorphism class of extensions, the bigger
will be the automorphism group of these extensions. When applying the reduction algorithm to a
polynomial of degree $p$ generating $L/K$, we have that when the unique ramification jump $t$ is an
integer and the additive polynomial $S_t(T)$ is not surjective we can do \emph{less simplifications}
to the coefficients of the Eisenstein polynomial.  Since any Eisenstein polynomial generating $L$ is
a possible output of the reduction algorithm (for a suitable choice of the $A_i$) we have that the
set of possible polynomials generating $L$ turns out to be ``smaller'', and that $L/K$ is Galois
having some non trivial automorphism and degree $p$.

For higher degree, and in particular when there are more ramification breaks, it becomes difficult
to generalize this observation, because a modification which appear to be trivial at the first order
may actually provoke some change to the higher order terms in the expansions. This fact also
justifies the claim that reducing the family to have exactly one polynomial for each isomorphism
class appears to be at least as hard as the computation of the number of isomorphisms for the
extension determined by one Eisenstein polynomial.

\section{A criterion to rule out possible reductions}
To complement the above reduction algorithm we give a synthetic criterion to exclude an Eisenstein
polynomial from generating an extension of which we know the set of all the reduced polynomials. In
particular given two polynomials $f(T)$ and $g(T)$ we can often rule out early the possibility that
a sequence of substitutions $\pi\rightarrow\pi+\theta\pi^{m+1}+\cO(\pi^{m+2})$, starting from level
$m=r$ say, may transform the minimal polynomial $f(T)$ of $\pi$ into the new minimal polynomial
$g(T)$, without having to compute the complete reduction.

Let's consider the monomial $(f_i-g_i)\pi^i$ having smallest valuation, which determines the
valuation of $f(\pi)-g(\pi)$, and assume its valuation to be equal to $v=n(\phi_{L/K}(r)+1)$ for
some real number $r$. Let's select sets of representatives $A_m$ which make $g(T)$ reduced, then we
say that $f(T)$ can be reduced to $g(T)$ \emph{greedily} if $g(T)$ is a possible output of the
reduction algorithm applied to $f(T)$ starting from step $m=r$. The proof of the following
proposition is clear.

\begin{prop}
\label{sec1_prop4}
  If $r$ is not an integer than $f(T)$ cannot be reduced greedily to $g(T)$.
\end{prop}

We also have the following Proposition, whose proof is immediate as well.

\begin{prop}
\label{sec1_prop5}
  If $r$ is an integer equal to a lower ramification break and $\overline{(f_i-g_i)\pi^{i-v}}$ is
  not in the image of $S_{t_i}(T)$, then $f(T)$ cannot be reduced greedily to $g(T)$.
\end{prop}

These observations are well complemented by the following Proposition, which makes them particularly
effective in the case of Galois extensions.

\begin{prop}
\label{sec1_prop6}
  Assume that one of $f(T)$ or $g(T)$ is known to generate a Galois extension, then $f(T)$ and
  $g(T)$ generate the same extensions if and only if one polynomial can be greedily reduced to the
  other.
\end{prop}

\begin{proof}
  For a suitable choice of representatives $g(T)$ is already reduced, and for Galois extensions
  there is only one reduced polynomial in view of Theorem \ref{theo1}, so applying greedily the
  reduction algorithm to $f(T)$ we obtain $g(T)$ as unique possible output. The other implication is
  clear.
\end{proof}

In other words for Galois extensions if $g(T)$ can be obtained in some way from $f(T)$, then it can
also be obtained in the greedy way.

When considering Galois extensions over $\Q_p$ the $S_{t_i}(T)$ are the zero map over the residue
field $\F_p$, so if $r$ is a ramification break then $f(T)$ and $g(T)$ certainly generate
non-isomorphic extensions, and we essentially recovered the main result of
\cite{yoshida2011ultrametric}.

However, it is possible to give a deeper criterion, which is more selective than what it is possible
by an inspection of $f(\pi)-g(\pi)$ at the first order.

Consider the range of monomials $f_{i,j}\pi_K^jT^i$ corresponding to one ramification break as
described in Prop. \ref{prop_range}, then the intuitive idea is that if we can obtain $f(T)$ from
$g(T)$ applying reductions of parameter $m\geq{}r$ then the first $r$ terms in each such interval
must be equal, because such terms are not going to be changed by any reduction of order
$\geq{}r$. Such ranges can be independently ``brought up to the front'' (with respect to the
$p$-adic valuation) computing formally a ramification polynomial of $f(T)-g(T)$, and considering the
coefficients of $T,T^p,T^{p^2},\dots$, as we can see observing the contributions to the coefficient
of $T^{p^\ell}$ in the ramification polynomial. Consequently taking into account a ramification
polynomial for $f(T)-g(T)$ provides a synthetic and effective formalism to describe how some sets of
coefficients must be equal in order to be able to pass from $f(T)$ to $g(T)$ via reduction step.

What we are going to prove is closely related to what was done in \cite{heiermann1996nouveaux} and
Theorem 4.6 in particular, and shares the philosophy that the sets of monomials $f_iT^i$ with a
fixed valuation of $i$ live an independent life from the other monomials, up to a certain extent,
and that when a uniformizer is changed $\pi\rightarrow\pi+\theta\pi^{m+1}+\cO(\pi^{m+2})$ the change
induced on minimal polynomial satisfies a certain continuity (in \cite{heiermann1996nouveaux} a
different kind of defining equation formed by a power series with coefficients in a set of
representatives was used rather than Eisenstein polynomials, but the underlying principle is the
same). From a more effective point of view, such a continuity provides an easily verifiable
criterion to exclude a polynomial from generating one fixed extension, which is particularly
effective in the case of Galois extensions thanks to Prop. \ref{sec1_prop6}. What we need seems not
to follow directly from the results of \cite{heiermann1996nouveaux} and additional steps would be
needed to switch to power series and back to Eisenstein polynomials, so we will avoid using the
slightly cumbersome notation of \cite{heiermann1996nouveaux} and prove our result directly.

For integers $a\geq0$ and $w$ let's define $P_{p^a}$, resp. $P_{p^a}(w)$, as the module generated
over $\cO_K$ by the monomials $cT^i$ such that $p^a\mid{}i$, resp. those monomials such that
additionally $v_L(c)+i\geq{}w$. If $cT^i\in{}P_{p^a}(w)$ for some $w$, than we have
\begin{equation}
\label{eq_deform}
   c\left(T+\theta{}T^{m+1}+\cO(T^{m+2})\right)^i - cT^i \in \sum_{j=0}^a P_{p^j}(w+e_L(a-j)+p^jm)
\end{equation}
as we can verify at once expanding the left hand side. Furthermore if $g\in{}P_{p^a}(w)$ and
$h\in{}P_{p^b}(z)$ than clearly we have $gh\in{}P_{p^{\min\{a,b\}}}(w+z)$.

Lets consider the ramification polygon $\Phi(T)=\pi^{-n}f(\pi{}T+\pi)$, then the coefficient of
$T^{p^a}$ is has valuation at least $\xi_{a}$. Assume $p^a\|i$, from a monomial $f_iT^i$ we have a contribution
$\binom{i}{p^a}f_i\pi^{i-n}T^{p^a}$ to the coefficient of $T^{p^a}$ in $\Phi(T)$, so $v_L(f_i)+i-n$
should be at least $\xi_a$, and consequently the monomial $f_iT^i$ is
contained in $P_{p^a}(\xi_a+n)$, being $v_L(f_i)+i\geq{}\xi_a+n$.

Consequently we have obtained that
\begin{equation}
\label{eq_structure}
 f(T) \in \sum_{j=0}^s P_{p^j}(\xi_j+n),
\end{equation}
where $s$ is the biggest integer such that $p^s\mid{}n$ (recall that $\xi_s=0$).

What observed above we obtain the following.

\begin{prop}
Let $F(T)=T+\theta_{m+1}T^{m+1}+\theta_{m+2}T^{m+2}+\dots$, then
\begin{equation}\label{eqcong}
   f(T) \equiv f(F(T)) \mod \sum_{j=0}^s P_{p^j}(\xi_j+n+p^j{}m).
 \end{equation}
\end{prop}

\begin{proof}
  Let's consider $f(T)-f(F(T))$, we will show that a monomial $f_iT^i$, which is contained in
  $P_{p^a}(\xi_a+n)$ by \eqref{eq_structure} say, yield various terms each having valuation at least
  $\xi_j+n+p^jm$ and in $P_{p^j}$, for some $j<a$. But we obtain terms in
  $P_{p^j}(\xi_a+n+e_L(a-j)+p^jm)$ by \eqref{eq_deform}, and $\xi_a+e_L(a-j)\geq\xi_j$ by Lemma
  \ref{sec1_lem2}.
\end{proof}

Assume $\pi=F(\rho)$ for a root $\rho$ of $g(T)$, we have now obtained a congruence property for the
power series $f(F(T))$, which clearly satisfies $f(F(\rho))=0$. The minimal polynomial of $\rho$ is
clearly a factor of $f(F(T))$ of degree $n$, and observe that the valuation the coefficient of $T^n$
is $0$ while the constant term has valuation $1$, so its Newton polygon has exactly one side of
length $n$ and slope $-\sfr{1}{n}$. In particular $g(T)$ is obtained by the factorization along the
Newton polygon, or equivalently collecting the roots with positive valuation, which is exactly what
is provided by the $p$-adic Weierstrass preparation Theorem.

We will however show a reduction which allows to approximate the minimal monic polynomial of $\rho$
starting from $f(F(T))$, and keeping the congruence \eqref{eqcong}. Let's start putting
$h_1(T)=f(F(T))$, and consider the polynomial $H_1(T)$ obtained taking the monomials of degree
$\geq{}n$ of $h_1(T)-T^n$.  If $H(T)=0$ then there is no such monomial, and $g(T)$ is a monic
polynomial of degree $n$, which is Eisenstein being $F(\pi)$ a root.

Let $cT^r$ a monomial of $H_1(T)$ which minimizes the quantity $v_L(c)+r$, and take the monomial
with $r$ as big as possible among those achieving the minimum of $v_L(c)+r$, which are in a finite
number. In other words, we are considering the higher valuation $\digamma$ on $\cO_L[[T]]$ where
\[
 \digamma(cT^r) = (\digamma_1(cT^r), \digamma_2(cT^r)) = (v_L(c)+r,\ -r)\in\Z^2
\]
and the elements of $\Z^2$ are ordered lexicographically, and we take $cT^r$ to be the monomial of
$H_1(T)$ which minimizes $\digamma(cT^r)$.

Let's replace now $h_1(T)$ with the new polynomial
\[
  h_2 = h_1(T) - h_1(T)\cdot cT^{r-n} = h_1(T)\cdot\left(1 - cT^{r-n}\right).
\]

We can see that applying iteratively such step $i$ times that either the minimum of the quantity
$v_L(c)+r$ for the monomials of degree $\geq{}n$ of $H_i(T)$ is increased, either is decreased the
biggest degree of the monomials achieving the minimum. Since the whole computation is done in
$\cO_K[[T]]$ the latter can only happen a finite number of times, and such minimum is increased
after a finite number of steps.

\begin{algorithm}
\caption{Lifting step}
\label{factor_lifting}
\begin{algorithmic}
\STATE $H_i(T)\gets\left\{\text{sum of monomials of degree $\geq{}n$ of }h_{i}(T)-T^n\right\}$
\STATE $cT^r\gets\left(\text{monomial of }H_i(T)\text{ minimizing }\digamma\right)$
\RETURN $h_i(T)\cdot(1-cT^{r-n})$
\end{algorithmic}
\end{algorithm}

After a sufficient number of iterations we can replace $h_i(T)$ with the polynomial $h(T)$ formed by
$T^n$ plus the monomials of degree $<n$ of $h_i(T)$, obtaining an Eisenstein polynomial such that
$h(F(\pi))$ is arbitrarily small, so $h(T)$ is itself an arbitrarily good approximation of the
minimal polynomial of $F(\pi)$.

We need to show that while the above procedure approximating $g(T)$ is carried on the congruence
satisfied by $f(F(T))$ is preserved. Indeed, assume the congruence to be satisfied by $h_i(T)$ and
assume $h_{i+1}(T)=h_i(T)(1-cT^{r-n})$. Then $cT^{r}\in{}P_{p^\ell}(\xi_\ell+n+p^\ell{}m)$ for some
$\ell$, being $cT^r$ a monomial of $H_i(T)$ and in view of the congruence which we assume to be
satisfied by $f(T)$ and $h_i(T)$. Let $bT^s$ be a monomial of $h_i(T)$, then
$bT^s\in{}P_{p^k}(\xi_k+n)$ for some $k$ by equation \eqref{eq_structure} and by the congruence
satisfied by $h_i(T)$. If $k<\ell$ we have
\[
    bT^s\cdot{}cT^{r-n} \in P_{p^k}(\xi_k+n+\xi_\ell+n+p^\ell{}m-n)\subseteq P_{p^k}(\xi_k+n+p^km),
\]
while when $k\geq\ell$ we have
\[
    bT^s\cdot{}cT^{r-n} \in P_{p^\ell}(\xi_k+n+\xi_\ell+n+p^\ell{}m-n)\subseteq P_{p^\ell}(\xi_\ell+n+p^{\ell}m).
\]
We obtained that subtracting $h_{i}(T)\cdot cT^{r-n}$ from $h_i(T)$ preserves the congruence.
Considering also an analogue of a ramification polynomial for $f(T)-g(T)$ we have the following
Theorem.

\begin{theo}
\label{theo2}
Let $f(T)$ be an Eisenstein polynomial of degree $n$ and $\pi$ a root, if $g(T)$ is another
Eisenstein polynomial of degree $n$ having $\rho\in{}K(\pi)$ as root, and
$\pi=\rho+\theta\rho^{m+1}+\cO(\rho^{m+2})$ then we have that
\[
   f(T) \equiv g(T) \mod \sum_{j=0}^s P_{p^j}(\xi_j+n+p^j{}m),
\]
and the polynomial
\[
   f(\pi+\pi{}T)-f(\pi)-g(\pi+\pi{}T)+g(\pi)
\]
has its Newton polygon contained in the Newton polygon of $f(\pi+\pi^{m+1}T)$.
\end{theo}

\begin{proof}
  We only need to prove the second assertion, but if $cT^r$ is in $P_{p^j}(\xi_j+n+p^j{}m)$ then for
  each $k\leq{}j$ the contribution of $c(\pi+\pi{}T)^r$ to the coefficient of $T^{p^k}$ has
  valuation at least $\xi_j+n+p^j{}m+(j-k)e_L$, which is at least $\xi_k+n+p^k{}m$ as shown above.
\end{proof}

\begin{remark}
  \label{rmk_resultant}
  We point out that the algorithm used during the proof to recover (an approximation of) the minimal
  polynomial of $\rho=F^{-1}(\pi)$ can be used to produce the minimal polynomial of a uniformizing
  element obtained deforming $\pi$ in a much quicker way than by computing a resultant
  $Res_U(g(U),T-(U+\theta{}U^{m+1}))$ as the determinant of a $(n+m)\times{}(n+m)$ matrix with
  coefficients in $\cO(T)$. Consequently taking $F(T)=T-\theta{}T^{m+1}$ and computing via the above
  approximation the minimal polynomial of the uniformizer $\rho$ such that
  $\pi=\rho-\theta\rho^{m+1}$, we obtain the minimal polynomial of a uniformizer
  $\rho=\pi+\theta\pi^{m+1}+\cO(\pi^{m+2})$, and this observation allows to exploit the free choice
  of $F(T)$ in Algorithm \ref{reduction_m} to avoid the computation of the resultant.
\end{remark}

\section{Construction of totally ramified class fields}

In this section we show how it is possible to convert a norm subgroup, representing a totally
ramified abelian extension via local class field theory, into the unique reduced Eisenstein
polynomial generating the extension.

We suppose given a finite index closed subgroup $N\subset{}K^\times$ such that $NU_{0,K}=K^\times$,
so that the corresponding extension by local class field theory is totally ramified. Being closed we
have $N\supset{}U_u$ for $u$ sufficiently big, this hypothesis is automatically satisfied when $K$
is a finite extension of $\Q_p$ and $N$ has finite index.

We assume $N$ to be described by a set of linear maps, one for each upper ramification break. That
is for all $u\geq0$ such that $U_{u,K}\nsubseteq{}NU_{u+1,K}$ we assume given a surjective
homomorphism
\[
        \nu_u : NU_{u,K} \rightarrow V_u
\]
having kernel exactly equal to $NU_{u+1,K}$, for some abstract group $V_u$ which is naturally an
$\F_p$-vector space for $u\geq1$. Take $\nu_u$ to be the trivial map to the trivial group $1$ when
$u$ is not an upper break, that is $NU_{u,K}=NU_{u+1,K}$. Note that the knowledge of all the maps
$\nu_u$ determines uniquely the group $N$.

The map $\nu_0$, when non-trivial, gives a condition on the representative $f_{0,1}$, or
equivalently on the residue class $\overline{f_0/\pi_K}$, this correspond to the well known explicit
description of local class field theory for tamely ramified extensions.  On the other hand the terms
appearing in a reduced polynomial in connection to the cokernels of the polynomials $S_{t_i}(T)$
attached to the lower breaks $t_i$ are all of the form $f_{0,j}\pi_K^j$, because the upper breaks
are integers by Hasse-Arf Theorem.

Consequently the choice of such representatives $f_{0,j}$ is determined by the condition that $f_0$
should be a norm from the extension determined by $N$, and a suitable $f_0$ can be selected changing
appropriately $\pi_K$.

The ramification data is described by the upper breaks $u\geq1$ and the dimensions of the
corresponding $V_u$. Consequently after selecting $f_0$ we have fixed a skeleton for the reduced
Eisenstein polynomial, formed by a well defined set of terms $f_{i,j}\pi_K^jT^i$, for $i\neq0$,
where the $f_{i,j}$ will be considered as unknowns in the set of representatives $R$. We will
describe how it is possible to recover the $f_{i,j}$ from the maps $\nu_u$.

The terms $f_{i,j}$ in a fixed range as in Prop.  \ref{prop_range} can be evaluated at the level
$\sigma_\ell$ when $p^\ell\|i$ say, making use the map $\nu_{\sigma_\ell/n-1}$ as we will now
show. If $m$ is such that $nj+i+p^\ell{}m=n\sigma_\ell+n$ it will be possible to describe the
dependence of $N_{K(\pi)/K}(1-\theta\pi^{m})$ on the coefficient $f_{i,j}$ at the first order,
obtaining a linear system from $\nu_{\sigma_\ell/n-1}$.

\begin{definition}
If $p^{\ell+1}|n$, we define $R_\ell$ to be the set of pairs $(i,j)$ such that $j\geq{}1$,
$0\leq{}i<n$, $p^\ell\|i$ and
\[
  \xi_\ell+n \leq nj+i < \sigma_\ell+n.
\]
We assume $R_\ell$ to be ordered depending on the value of $nj+i$. We define $M(i,j)$ to be the
number $m$ such that
\[
    nj+i +p^\ell{}m = \sigma_\ell+n.
\]
\end{definition}

We remark that if the extension is abelian then $n|\sigma_\ell$ by Hasse-Arf theorem, so the $m$
defined above is always an integer $\leq{}\tau_\ell$ and prime with $p$.

\subsection{Dependence of norms on a $f_{i,j}$}
We will now track the dependence of a norm $N_{K(\pi)/K}(1-\theta\pi^m)$ on a representative
$f_{i,j}$ appearing in the expansion of a coefficient. To do so, let's treat $f_{i,j}$ as an
indeterminate, and apply a sufficient number of steps of Algorithm \ref{factor_lifting} to pass from
$f(T+\theta{}T^{m+1})$ to the minimal polynomial $g(T)$ of $\rho$, where
$\pi=\rho+\theta\rho^{m+1}$.

We clearly have $\rho=\pi-\theta\pi^{m+1}+\cO(\pi^{m+2})$, and $g_0/f_0$ will be the norm of some
element which is $1-\theta\pi^{m}+\cO(\pi^{m+1})$. Since we will obtain that the changes induced changing $f_{i,j}$
on all $N_{K(\pi)/K}(1-\theta\pi^r+\cO(\pi^{r+1}))$ for $r>m$ will be even smaller $p$-adically,
this makes possible to ignore the $\cO(\pi^{m+1})$.

In the expansion of $f(T+\theta{}T^{m+1})$ a term $f_{i,j}\pi_K^j(T+\theta{}T^{m+1})^i$ appears, and
it has $f_{i,j}\pi_K^jT^i$ as main term. In the algorithm we start with
$h_0(T)=f(T+\theta{}T^{m+1})$, and at the $i$-th step we subtract $h_i(T)\cdot{}cT^{r-n}$ from
$h_i(T)$ for a monomial $cT^{r-n}$ in $P_{p^k}(\xi_k+p^km)$ for some $k\leq{}s$. From the
monomial $f_{i,j}\pi_K^jT^i$ the other terms in $f_{i,j}$ which may appear in the algorithm have
coefficient with valuation at least
\[
    nj+i +\min_{0\leq{}k\leq{}s}\{ \xi_k+p^km \} = nj + i + n\phi_{L/K}(m),
\]
where we are considering the mixed valuation $\digamma_1(\pi_K^aT^b)=na+b$ on $\cO_L[[T]]$. 

Note that the minimum of $\xi_k+p^k{}m$ is obtained as the minimum of the piecewise linear function
$nN(x)+mx$, which is $n\phi_{L/K}(m)$ in view of what proved before Proposition \ref{sec1_prop1}. We
will denote for convenience this quantity as
\[
   A_{i,j}(m) = nj+i+n\phi_{L/K}(m).
\]

The term $f_{i,j}\pi_K^jT^i$ is the main term coming from
$f_{i,j}\pi_K^j(T+\theta{}T^{m+1})^i$, and the second contribution can be found considering the expansion
\[
   (1+\theta{}T^m)^i = 1 + \binom{i}{p^\ell}(\theta{}T^m)^{p^\ell}
      + \binom{i}{p^{\ell-1}}(\theta{}T^m)^{p^{\ell-1}}+\dots.
\]
Putting as usual $p^\ell\|i$, we denote the valuation of the second term as
\[
   B_{i,j}(m) = nj+i+\min_{0\leq{}k\leq\ell}\left\{e_L(\ell-k) + p^km\right\}.
\]
Such term is equal to
\[
   \binom{i}{p^\ell}\pi_K^jT^i\cdot(\theta{}T^m)^{p^\ell}
\] 
as long as $mp^\ell<e_L+mp^{\ell-1}$, that is $m<\sfr{e_L}{(p^\ell-p^{\ell-1})}$, and if
$m\leq{}M(i,j)\leq\tau_\ell$ this condition is certainly satisfied because
\[
   m \leq \tau_\ell \leq \sfr{e_L}{(p^{\ell+1}-p^{\ell})}
\]
by Lemma \ref{sec1_lem2}.

By Proposition \ref{prop_range} we can assume $\ell<s$, and if $m\leq{}\tau_\ell$ we always have
\[
  B_{i,j}(m) = nj+i+p^\ell{}m < A_{i,j}(m).
\]
So the main contribution to $N_{K(\pi)/K}(1-\theta\pi^m)=g_0/f_0$ originated from
$f_{i,j}\pi_K^jT^i$ is only coming from $f_{i,j}\binom{i}{p^\ell}\pi_K^jT^i\cdot(\theta{}T^m)^{p^\ell}$.

When $m=M(i,j)$ we obtain a condition on $f_{i,j}$ to have $N_{K(\pi)/K}(1-\theta\pi^m)$ in
the norm group for each $\theta$, which can be used to determine the representative $f_{i,j}$.

This can be made to work when only one representative $f_{i,j}$ is unknown, but a more refined study
is needed if we have to determine them all. In particular, we will see that there exists an ordering
of such unknowns that allows to determine them all inductively. What complicates this idea is that
it turns out to be necessary to interleave in a suitable way the ranges of representatives
considered in Prop. \ref{prop_range}.

It will be convenient to write down a comfortable lower bound for the functions $A_{i,j}$ and
$B_{i,j}$ to ensure that the value of a particular representative $f_{i,j}$ has no influence on
$N_{K(\pi)/K}(\pi-\theta\pi^{m+1})$ modulo a suitable power of $\pi$. In particular we have that
they are all bounded by
\[
  C_{i,j}(m) = nj+i + \min_{0\leq{}k\leq{}\ell}\left\{\xi_k-\xi_\ell+p^km\right\},
\]
where $\ell=v_p(i)$ as usual. We resume the properties proved in the following Lemma.

\begin{lemma}
\label{th3lem1}
  Denote with $\pi$ a root of $f(T)$, and let $(i,j)\in{}R_\ell$. Then,
  for each $\theta\in{}U_K$, the value of
\[
   N_{K(\pi)/K}(1-\theta\pi^m) \mod \kp_K^u
\]
does not depend on $f_{i,j}$, whenever $u$ is $\leq{}C_{i,j}(m)/n-1$. Furthermore if $m=M(i,j)$ then
$C_{i,j}(m)=\sigma_\ell+n$, and putting $u=C_{i,j}(m)/n-1=\phi_{L/K}(\tau_\ell)$ we have
\begin{equation}
\label{normdelta}
   N_{K(\pi)/K}(1-\theta\pi^m) = N_{K(\pi_0)/K}(1-\theta\pi_0^m)
     + \pi_K^uf_{i,j}\lambda_{i,j}\theta^{p^\ell} + \cO(\pi_K^{u+1}),
\end{equation}
where $\lambda_{i,j}$ is a fixed unit defined as
\[
   \lambda_{i,j} = \binom{i}{p^\ell}\cdot(\sfr{-f_0}{\pi_K})^{(i+p^\ell{}m)/n-1},
\]
and we set $\pi_0$ to be a root of the polynomial obtained from $f(T)$ setting $f_{i,j}$ to $0$.
\end{lemma}

\begin{proof}
  We just have to prove the \eqref{normdelta}. For $m=M(i,j)$ the variation of the constant term
  comes from the monomial $f_{i,j}\binom{i}{p^\ell}\pi_K^jT^i\cdot(\theta{}T^m)^{p^\ell}$ in the
  expansion of $f(T+\theta{}T^{m+1})$, and during the reduction each $T^n$ is transformed into
  $-f_0$. Dividing by $f_0$ we obtain that the variation for $N_{K(\pi)/K}(1-\theta\pi^m)$, which
  modulo $\pi_K^{u+1}$ is
\[
   f_{i,j}\pi_K^j\binom{i}{p^\ell}(-f_0)^{(i+p^\ell{}m)/n-1}\theta^{p^\ell} 
     = \pi_K^uf_{i,j}\lambda_{i,j}\theta^{p^\ell}. \qedhere
 \]
\end{proof}

We will now assume that some $f_{i,j}$ have been determined and some not yet, and will show that it
is possible to determined some of the unknown ones. For $\ell$ such that $p^{\ell+1}\|n$, consider
the range of terms $R_\ell$ like in Prop. \ref{prop_range}, and let $(i_\ell,j_\ell)$ be the
smallest pair $(i,j)\in{}R_\ell$ (i.e. the pair in $R_\ell$ with $nj+i$ as small as possible) such
that the corresponding $f_{i,j}$ has not been identified yet.

We first prove a couple of technical lemmas about the functions $C_{i,j}(x)$.

\begin{lemma}
\label{cijtech1}
The functions $C_{i,j}(x)$ are strictly increasing, and for $(i,j)\neq(i',j')$ then the functions
$C_{i,j}(x)$ and $C_{i',j'}(x)$ are always different except possibly at one point. The difference
$C_{i,j}(x)-C_{i',j'}(x)$ is constant when $v_p(i)=v_p(i')$, and $C_{i',j'}(x)$ can surpass
$C_{i,j}(x)$ only when $v_p(i')>v_p(i)$.
\end{lemma}
\begin{proof}
Follows directly from the definition.
\end{proof}

\begin{lemma}
\label{cijtech2}
Let $\Pi_L(x)$ be the real function $x\mapsto\min\{px,x+e_L\}$, and let $\Pi^{[h]}$ denote the
$h$-times composition. Then if $y=C_{i,j}(m)-n$, for $(i,j)\in{}R_\ell$ and some $m\in{}\mathbb{N}$,
then we have
\[
   \Pi^{[h]}(y) \geq C_{i,j}(p^hm)-n,
\]
for each $h\geq1$.
\end{lemma}
\begin{proof}
  It's enough to prove that $\Pi(y)\geq{}C_{i,j}(pm)-n$, and we can do so proving that both $py$ and
  $y+e_L$ are bigger. This follows easily from the definition and from Lemma \ref{sec1_lem2}.
\end{proof}

Let $K(x)$ be the function defined as
\begin{equation}
\label{mincij}
   K(x) = \min_{0\leq\ell<s} C_{i_\ell,j_\ell}(x),
\end{equation}
it is again a strictly increasing function which for each $m$ describes up to which precision we can
compute $N_{K(\pi)/K}(\pi-\theta\pi^{m+1})$, with the given information about the coefficients of
$f(T)$.

Let $I_\ell$ be the set of real $x$ where $K(x)=C_{i_\ell,j_\ell}(x)$, since two functions
$C_{i,j}(x)$ can only cross once we have that $I_\ell$ is a (possibly infinite) topologically closed
real interval, and taking into account the conditions under which a surpass may happen of Lemma
\ref{cijtech1} we obtain that $I_k$ lies before $I_\ell$ if $k>\ell$. Let's merge in a unique bigger
interval the $I_\ell$ such that the value of $\tau_\ell$ is the same, and let $J_r$ be formed by the
union of the $I_\ell$ such that $\tau_\ell=t_r$. Let's also put $k_r=K^{-1}(n\phi(t_r)+n)$ for each
$r$, then $k_r$ is contained in the \emph{interior} of $J_r$ for some $r$, thanks to the following
Lemma.

\begin{lemma}
  Let $A_1,\dots,A_m$ be a sequence of intervals with extrema $\R\cup\{\pm\infty\}$, such that
  $A_{i+1}$ begins exactly where $A_i$ ends. Let $a_1<a_2<\dots<a_m$ be real numbers contained in
  the interior of $\bigcup_{i=1}^m A_i$. Then $a_i$ is contained in the interior of $A_i$, for some
  $1\leq{}i\leq{}m$.
\end{lemma}
\begin{proof}
  The thesis is trivial when $m=1$. If $m>1$, then either $a_m$ is in the interior of $A_m$, either
  we have that $a_1,a_2,\dots,a_{m-1}$ are contained in the interior of $\bigcup_{i=1}^{m-1}A_i$ and
  the thesis follows by induction.
\end{proof}

So, let $r$ be such that $k_r$ is in the interior of $J_r$, and let $L_1\leq{}L_2$ be the integers
such that $\tau_\ell=t_r$ if and only if $\ell\in[L_1,L_2]$. Then we have that
$C_{i_\ell,j_\ell}(k_r)$ is at least $n\phi_{L/K}(t_r)+n$ for $\ell\in[L_1,L_2]$, and strictly bigger
for $\ell\notin[L_1,L_2]$. Put $m=k_r$.

For $L_1\leq\ell\leq{}L_2$ let's redefine $(i_\ell,j_\ell)$ to be the unique pair $(i,j)$ in
$R_\ell$ that $C_{i,j}(m)=n\phi_{L/K}(t_r)+n$. Then some pairs $(i_\ell,j_\ell)$ will be unchanged
while some others may be set to correspond to representatives $f_{i_\ell,j_\ell}$ which are already
known. This makes no harm since all these representatives will now be determined
simultaneously. Since the equality need to hold already for some of the original $(i_\ell,j_\ell)$
we obtain that $m=M(i_\ell,j_\ell)$ is an integer $<\tau_\ell=t_r$ and prime with $p$.

Let's consider the terms $f_{i_\ell,j_\ell}\pi_K^{j_\ell}T^{i_\ell}$ for $\ell$ in the given range,
and lets vary the $f_{i_\ell,j_\ell}$. Put $u=\phi_{L/K}(t_r)$. 

\begin{lemma}
If the $f_{i,j}$ already known were determined inductively we already have
\[
   N_{K(\pi)/K}(1-\theta\pi^{m}) \in U_uN,
\]
and the class modulo $U_{u+1}N$ only depend on $\bar\theta$.
\end{lemma}
\begin{proof}
  Indeed let $u'=\phi_{L/K}(t_{r'})\leq{}u$ be an upper ramification jump and $m'\geq{}m$, and
  assume at least one of these inequality to be strict. We will show that there was a previous step
  where we computed some currently known coefficient, by requesting $N_{K(\pi)/K}(1-\theta\pi^{m'})$
  to be in $U_{u'+1}N$ for each $\theta$.

If $\phi_{L/K}(m')>u'$ then we always have $N_{K(\pi)/K}(1-\theta\pi^{m'})\in{}U_{u'+1}$, by the
properties of the norm map \cite[Chap 3, Prop 3.1]{fesenko2002local}. Consequently assume
$\phi_{L/K}(m')\leq{}u'$, or $m'\leq{}\psi_{L/K}(u')=t_{r'}$ applying $\psi_{L/K}$. If $(m,p)=1$,
let $[L_1',L_2']$ be the interval of possible $\ell$ such that $\tau_\ell=t_{r'}\geq{}m$, and for
each $\ell$ we have a pair $(i_\ell',j_\ell')$ for $L_1'\leq\ell\leq{}L_2'$ which has the property
that $m'=M(i_\ell',j_\ell')$, and $C_{i_\ell',j_\ell'}(m')=nu'+n\leq{}K(m)$. Thus the
$C_{i_\ell',j_\ell'}(x)$ is certainly not $\geq{}K(x)$ and cannot appear in the \eqref{mincij}. This
means that the $f_{i_\ell',j_\ell'}$ for $L_1'\leq\ell\leq{}L_2'$ have been determined at a previous
step, where we guaranteed that $N_{K(\pi)/K}(1-\theta\pi^{m'})\in{}U_{u'+1}$.

When $p|m$ let's consider the elements of the form $(1+\theta\pi^{m''/p^w})^{p^w}$ as generators of
$U_{m'}/U_{m'+1}$ (see \cite[Chap. 1, Prop. 5.7]{fesenko2002local}, and note that
$m'\lneq\sfr{pe_L}{(p-1)}$). We are done if we show that $N(1+\theta\pi^{m''})\in{}NU_{u''+1}$ for
$m''=m'/p^w$ and each upper break $u''\leq\Pi_K^{[-w]}(u')$, where $\Pi_K(x)=\min\{px,x+e_K\}$. 

We can assume $m''\leq{}\psi_{L/K}(u'')=t_{r''}$ as above. Consider as above an $\ell$ such that
$\tau_\ell=t_{r''}=\psi_{L/K}(u'')$, and a pair $(i_{\ell}'',j_\ell'')\in{}R_\ell$ with
$m''=M(i_{\ell}'',j_\ell'')$. Then $C_{i_\ell'',j_\ell''}(m'')=nu''+n$, and
$C_{i_\ell'',j_\ell''}(m')$ is certainly $\leq{}\Pi_L^{[w]}(nu'')+n\leq{}nu'+n$ by Lemma
\ref{cijtech2}. We have that $C_{i_\ell'',j_\ell''}(x)$ is not $\geq{}K(x)$, and the condition
$N(1+\theta\pi^{m''})\in{}NU_{u''+1}$ was verified in a previous step.
\end{proof}

We have from the Lemma \ref{th3lem1}, applied once for each pair $(i_\ell,j_\ell)$ for
$L_1\leq\ell\leq{}L_2$, that changing the $f_{i_\ell,j_\ell}$ we have
\begin{align*}
   N_{K(\pi)/K}(1-\theta\pi^{m}) = N_{K(\pi_0)/K}(1-\theta\pi_0^{m}) +
    \sum_{\ell=L_1}^{L_2}
    \pi_K^uf_{i_\ell,j_\ell}\lambda_{i_\ell,j_\ell}\theta^{p^\ell}
   + \cO\left(\pi^{u+1}\right),
\end{align*}
where $\pi_0$ is the root of the polynomial with the unknown $f_{i_\ell,j_\ell}$ all set to $0$.

Now for some choice of the $f_{i_\ell,j_\ell}$ we want $N_{K(\pi)/K}(1-\theta\pi^m)$ be in
the kernel of $\nu_u:NU_u/NU_{u+1}\stackrel{\sim}{\longrightarrow}{}V_u$, for each residue
representative $\theta$. The condition only depends on the reductions $\bar\theta$ and
$\overline{f_{i_\ell,j_\ell}}$ and is $\F_p$-linear in them, so we can impose it to hold only for a
set of values for $\bar\theta$ which generate $\kk_K$ over $\F_p$. Let
$\bar\alpha_1,\dots,\bar\alpha_F$ be such a basis, then we can decompose each
\[
   \overline{f_{i_\ell,j_\ell}}=\sum_{k=1}^F\overline{f_{i_\ell,j_\ell,k}}\bar\alpha_k
\]
and consider $\overline{f_{i_\ell,j_\ell,k}}$ as unknown over $\F_p$.

Fix a lifting $\alpha_k\in\cO_K$ for each $\bar\alpha_k$, we have an equation
\[
  \nu_u(N_{K(\pi_0)/K}(1-\alpha_k\pi_0^{m}) +
    \sum_{\ell=L_1}^{L_2} \pi_K^u\left[
    \sum_{k=1}^f f_{i_\ell,j_\ell,k}\alpha_k\right]\lambda_{i_\ell,j_\ell}\alpha_k^{p^\ell}) = 0
\]
in $V_u$ for each generator $\bar\theta=\bar\alpha_k$. The codomain $V_u$ has $\F_p$-dimension equal
to the ``multiplicity'' of the ramification jump $\Lambda=L_2-L_1+1$. So we have an inhomogeneous
system formed by $F\cdot\Lambda$ equations over $\F_p$, for the same number of unknowns
$\overline{f_{i_\ell,j_\ell,k}}$.

Consequently we have a unique solution for the $\overline{f_{i_\ell,j_\ell,k}}$ if we can prove that
the system is non-degenerate. We can equivalently prove that the connected homogeneous system
\[
  \nu_u\left(1 +
    \sum_{\ell=L_1}^{L_2} \pi_K^u\left[\sum_{k=1}^f
      f_{i_\ell,j_\ell,k}\alpha_k\right]\lambda_{i_\ell,j_\ell}\theta^{p^\ell}\right) = 0,
  \qquad\text{for all }\theta
\]
has no non-trivial solution. Subtracting $1$ and dividing by a suitable power of $\pi_K$, the system
can be interpreted as the request that the additive polynomial
\[
   \sum_{\ell=L_1}^{L_2} \overline{f_{i_\ell,j_\ell}}\overline{\lambda_{i_\ell,j_\ell}}\bar\theta^{p^\ell}
\]
in $\bar\theta$ should have range identically contained in a subspace of codimension
$L_2-L_1+1$. Since the powers of $\bar\theta$ appearing are $p$-th powers ranging from $p^{L_1}$ to
$p^{L_2}$, we have that its corank as linear map is at most $L_2-L_1$ (see \cite[Chap. 5,
\S2]{fesenko2002local}). Consequently such a non-trivial solution of the homogeneous system is
impossible, and the original system is non-degenerate.

\begin{theo}
  Given a closed finite index subgroup $N\subset{}K^{\times}$ corresponding to a totally ramified
  class field, there exists an ordering of the representatives $f_{i,j}$ appearing in the expansion
  of the coefficients of a generic Eisenstein polynomial which allows to determine the coefficients
  of the reduced Eisenstein polynomial generating the extensions corresponding to $N$.
\end{theo}

It is indeed clear that the above procedure where the $f_{i,j}$ are obtained solving linear
equations can be converted into an algorithm to construct explicitly a minimal equation
corresponding to a class field. The $f_{i,j}$ allowed by the ramification data but with
$(i,j)\notin{}R_\ell$ for each $\ell$ can be assumed to be all $0$, or set to any arbitrary value
(indeed, they are exactly those set by the reduction algorithm). On the other hand, during the
construction we guarantee that for each $m>0$ and prime with $p$, and for each $\theta$, we have
$N_{K(\pi)/K}(1-\theta\pi^m)\in{}N$, because the condition $\nu_u(N_{K(\pi)/K}(1-\theta\pi^m))=0$ is
verified for all $u$ at some point of the algorithm, implying that all norms are contained in $N$.

\begin{remark}
  We observe that this construction produces an alternative and constructive proof of the Existence
  Theorem of class field theory for totally ramified extensions, because for each finite index
  closed subgroup $N$ of $K^{\times}$ with index $n$ we construct an extensions of degree $n$ having
  norm subgroup contained in $N$. We can construct precisely one reduced polynomial of degree
  $(K^\times:N)$ for each $N$, and considering all possible reduction steps we have easily that the
  group of norms has to be exactly equal to $N$, and that all intermediate fields
  $K_{t_r^+}/K_{t_r}$ are Galois, so the generated field $L/K$ has to be Galois by Theorem
  \ref{theo1}.

  It would be interesting to extend this construction to recover the Artin map from $K^\times/N$ to
  $\Gal(K/L)$, proving that these two groups are indeed isomorphic and describing explicitly the
  isomorphism. The above methods do not even give an easy proof that the extension obtained is
  abelian, without assuming local class field theory.
\end{remark}

\bibliographystyle{amsalpha}
\bibliography{biblio}

\end{document}

%% file: newton1.pdf_tex

\begingroup
  \makeatletter
  \providecommand\color[2][]{%
    \errmessage{(Inkscape) Color is used for the text in Inkscape, but the package 'color.sty' is not loaded}
    \renewcommand\color[2][]{}%
  }
  \providecommand\transparent[1]{%
    \errmessage{(Inkscape) Transparency is used (non-zero) for the text in Inkscape, but the package 'transparent.sty' is not loaded}
    \renewcommand\transparent[1]{}%
  }
  \providecommand\rotatebox[2]{#2}
  \ifx\svgwidth\undefined
    \setlength{\unitlength}{240pt}
  \else
    \setlength{\unitlength}{\svgwidth}
  \fi
  \global\let\svgwidth\undefined
  \makeatother
  \begin{picture}(1,0.66666667)%
    \put(0,0){\includegraphics[width=\unitlength]{newton1.pdf}}%
    \put(0.90666667,0.01666666){\color[rgb]{0,0,0}\makebox(0,0)[lb]{\smash{$\gamma_0$}}}%
    \put(0.44,0.01666666){\color[rgb]{0,0,0}\makebox(0,0)[lb]{\smash{$\gamma_1$}}}%
    \put(0.19333333,0.01666666){\color[rgb]{0,0,0}\makebox(0,0)[lb]{\smash{$\gamma_2$}}}%
    \put(0.04666667,0.01666666){\color[rgb]{0,0,0}\makebox(0,0)[lb]{\smash{$\gamma_3$}}}%
    \put(0.66,0.13333332){\color[rgb]{0,0,0}\makebox(0,0)[lb]{\smash{$-\sfr{t_1}{n}$}}}%
    \put(0.3,0.29999999){\color[rgb]{0,0,0}\makebox(0,0)[lb]{\smash{$-\sfr{t_2}{n}$}}}%
    \put(0.15,0.49999999){\color[rgb]{0,0,0}\makebox(0,0)[lb]{\smash{$-\sfr{t_3}{n}$}}}%
  \end{picture}%
\endgroup

%% file: copolygon.pdf_tex

\begingroup
  \makeatletter
  \providecommand\color[2][]{%
    \errmessage{(Inkscape) Color is used for the text in Inkscape, but the package 'color.sty' is not loaded}
    \renewcommand\color[2][]{}%
  }
  \providecommand\transparent[1]{%
    \errmessage{(Inkscape) Transparency is used (non-zero) for the text in Inkscape, but the package 'transparent.sty' is not loaded}
    \renewcommand\transparent[1]{}%
  }
  \providecommand\rotatebox[2]{#2}
  \ifx\svgwidth\undefined
    \setlength{\unitlength}{296pt}
  \else
    \setlength{\unitlength}{\svgwidth}
  \fi
  \global\let\svgwidth\undefined
  \makeatother
  \begin{picture}(1,0.59459459)%
    \put(0,0){\includegraphics[width=\unitlength]{copolygon.pdf}}%
    \put(0.78378378,-0.00000001){\color[rgb]{0,0,0}\makebox(0,0)[lb]{\smash{$t_3 (\tau_0)$}}}%
    \put(0.64864865,-0.00000001){\color[rgb]{0,0,0}\makebox(0,0)[lb]{\smash{$t_2 (\tau_1)$}}}%
    \put(0.42702703,-0.00000001){\color[rgb]{0,0,0}\makebox(0,0)[lb]{\smash{$t_1 (\tau_2)$}}}%
    \put(0.31408056,0.25050428){\color[rgb]{0,0,0}\makebox(0,0)[lb]{\smash{$\scriptstyle{\sigma_1}$}}}%
    \put(0.31360016,0.46185401){\color[rgb]{0,0,0}\makebox(0,0)[lb]{\smash{$\scriptstyle{\sigma_2}$}}}%
    \put(0.26756757,0.52432432){\color[rgb]{0,0,0}\makebox(0,0)[lb]{\smash{$\scriptstyle{\sigma_3}$}}}%
    \put(0.18918919,0.48648648){\color[rgb]{0,0,0}\makebox(0,0)[lb]{\smash{\framebox{$p\nmid{}i$}}}}%
    \put(0.01081081,0.3108108){\color[rgb]{0,0,0}\makebox(0,0)[lb]{\smash{$f_{i,j}\pi_K^jT^i$}}}%
    \put(0.19459459,0.39459459){\color[rgb]{0,0,0}\makebox(0,0)[lb]{\smash{\framebox{$p\|{}i$}}}}%
    \put(0.2027027,0.21351351){\color[rgb]{0,0,0}\makebox(0,0)[lb]{\smash{\framebox{$p^2\|{}i$}}}}%
    \put(0.40000528,0.09498554){\color[rgb]{0,0,0}\makebox(0,0)[lb]{\smash{$p^3$}}}%
    \put(0.55405405,0.33783783){\color[rgb]{0,0,0}\makebox(0,0)[lb]{\smash{$p^2$}}}%
    \put(0.72162162,0.47297297){\color[rgb]{0,0,0}\makebox(0,0)[lb]{\smash{$p$}}}%
    \put(0.87837838,0.5081081){\color[rgb]{0,0,0}\makebox(0,0)[lb]{\smash{$1$}}}%
    \put(0.26756757,0.45405405){\color[rgb]{0,0,0}\makebox(0,0)[lb]{\smash{$\scriptstyle{\xi_3}$}}}%
    \put(0.31408056,0.16942319){\color[rgb]{0,0,0}\makebox(0,0)[lb]{\smash{$\scriptstyle{\xi_1}$}}}%
    \put(0.31408056,0.34780157){\color[rgb]{0,0,0}\makebox(0,0)[lb]{\smash{$\scriptstyle{\xi_2}$}}}%
  \end{picture}%
\endgroup

%% file: art5.bbl
\providecommand{\bysame}{\leavevmode\hbox to3em{\hrulefill}\thinspace}
\providecommand{\MR}{\relax\ifhmode\unskip\space\fi MR }
\providecommand{\MRhref}[2]{%
  \href{http://www.ams.org/mathscinet-getitem?mr=#1}{#2}
}
\providecommand{\href}[2]{#2}
\begin{thebibliography}{Yam68}

\bibitem[Ama71]{Amano1971}
Shigeru Amano, \emph{Eisenstein equations of degree {$p$} in a
  {$\mathfrak{p}$}-adic field}, J. Fac. Sci. Univ. Tokyo Sect. IA Math
  \textbf{18} (1971), 1--21.

\bibitem[Del84]{deligne1984lescorps}
Pierre Deligne, \emph{Les corps locaux de caract\'eristique {$p$}, limites de
  corps locaux de caract\'eristique {$0$}}, Representations of reductive groups
  over a local field, Travaux en Cours, Hermann, Paris, 1984, pp.~119--157.

\bibitem[FV02]{fesenko2002local}
Ivan~Borisovich Fesenko and Sergei~Vladimirovich Vostokov, \emph{{Local fields
  and their extensions}}, American Mathematical Society, 2002.

\bibitem[Hei96]{heiermann1996nouveaux}
Volker Heiermann, \emph{De nouveaux invariants num{\'e}riques pour les
  extensions totalement ramifi{\'e}es de corps locaux}, Journal of Number
  Theory \textbf{59} (1996), no.~1, 159--202.

\bibitem[Hel91]{Helou1991}
Charles Helou, \emph{{On the ramification breaks}}, Communications in Algebra
  \textbf{19} (1991), no.~8, 2267--2279.

\bibitem[Kra62]{Krasner1962}
Marc Krasner, \emph{{Nombre des extensions d'un degr\'{e} donn\'{e} d'un corps
  $\mathfrak{p}$-adique}}, C. R. Acad. Sc. Paris \textbf{254} (1962),
  3470--3472, \textit{Ibidem} 255:224--226, 1682--1684, 2342--2344, 3095--3097,
  1962.

\bibitem[Li97]{li1997p}
Hua-Chieh Li, \emph{p-adic power series which commute under composition},
  Transactions of the American Mathematical Society \textbf{349} (1997), no.~4,
  1437--1446.

\bibitem[Lub81]{lubin1981local}
Jonathan Lubin, \emph{The local {K}ronecker-{W}eber theorem}, Transactions of
  the American Mathematical Society \textbf{267} (1981), no.~1, 133--138.

\bibitem[Pau06]{pauli2006constructing}
Sebastian Pauli, \emph{Constructing class fields over local fields}, Journal de
  th{\'e}orie des nombres de Bordeaux \textbf{18} (2006), no.~3, 627--652.

\bibitem[PR01]{Pauli2001computation}
Sebastian Pauli and Xavier-Fran\c{c}ois Roblot, \emph{{On the computation of
  all extensions of a $p$-adic field of a given degree}}, Mathematics of
  Computation \textbf{70} (2001), no.~236, 1641--1660.

\bibitem[Ser78]{serre1978formule}
Jean-Pierre Serre, \emph{{Une ``formule de masse'' pour les extensions
  totalement ramifi{\'e}es de degr{\'e} donn{\'e} d'un corps local}}, C.R.
  Acad. Sci. Paris S\'er. A-B \textbf{286} (1978), 1031--1036.

\bibitem[Ser79]{serre1979local}
\bysame, \emph{{Local fields}}, Springer, 1979.

\bibitem[Yam68]{Yamamoto1968}
Sunao Yamamoto, \emph{{On a property of the Hasse's function in the
  ramification theory}}, Memoirs of the Faculty of Science, Kyushu University.
  Series A, Mathematics \textbf{22} (1968), no.~2, 96--109.

\bibitem[Yos11]{yoshida2011ultrametric}
Manabu Yoshida, \emph{An ultrametric space of {E}isenstein polynomials and
  ramification theory}, Arxiv preprint arXiv:1105.5221 (2011).

\end{thebibliography}
